\documentclass[11pt,reqno]{amsart}
\usepackage{amssymb,color, fullpage}

\title{Boole's Chapter XV: Syllogism Details}
\author{Stanley Burris}
\email{snburris at uwaterloo.ca}
\date{\today}                                           

\allowdisplaybreaks
\newcommand{\RR}[1]{\begin{color}{red} #1 \end{color}}

\begin{document}
\maketitle

\begin{abstract}
In Boole's famous 1854 book {\em The Laws of Thought\/} the mathematical analysis of Aristotelian 
logic was relegated to Chapter XV, the last chapter before his treatment of probability theory. 
This chapter is Boole's tour de force to show that he had a uniform method to obtain all valid
syllogisms in his version of Aristotelian logic, namely he applied {\em reduction}, {\em elimination}
and {\em solution} in that order to equational expressions for the premises.
The premises of a syllogism were expressed as a pair of equations in 7 variables, but then
all algebraic steps between this and the final expressions for $x$, $vx$ and $1-x$ were omitted.
The somewhat tedious details of those missing steps are given in this note. 
It is assumed that the reader is
familiar with Boole's reduction, elimination and solution theorems.
\end{abstract}

\section{Boole's two pairs of equations}

In LT \cite{Boole-1854} ({\em The Laws of Thought}), pp.~232--236, in order to justify the 
simple rules describing valid syllogisms that were stated in his 1848 paper \cite{Boole-1848},
Boole presented three solutions 
$$
\begin{array} {r c l}
x &=& f_1(v,v',w,w')z +  f_0(v,v',w,w')(1-z)\\
vx &=& g_1(v,v',w,w')z +  g_0(v,v',w,w')(1-z)\\
1-x &=& h_1(v,v',w,w')z +  h_0(v,v',w,w')(1-z)
\end{array}
$$
for each of the two pairs of equations
$$
\mathrm{I}.\ 
\left\{\begin{array}{r c l}
vx&=&v'y \\
wz&=&w'y
\end{array} \right.
\qquad \qquad
\mathrm{II}.\ 
\left\{\begin{array}{r c l}
vx&=&v'y\\
wz&=&w'(1-y) ,
\end{array} \right.
$$
after eliminating $y$.

By independently substituting $1-x$ for $x$, $1-y$ for $y$ and $1-z$ for $z$
in I and II, carrying the substitutions over to the three solutions, 
Boole said that one would have the equational expressions for all possible 
arguments with categorical premises. However each of the three conclusions obtained by
solution would not correspond to a categorical proposition unless either the
coefficient of $z$ or the coefficient of $1-z$ vanished.

In order for these equational arguments to be applicable to the determination
 of valid syllogisms he needed, in each of the three solutions 
in each of the two cases,
 to find all {\em permissible} substitutions of 1 for some of the variables $v,v',w,w'$ such that 
 the coefficient of $z$ or the coefficient of $1-z$ would be 0. By a permissible substitution
 is meant one where not both $v,v'$ nor both $w,w'$ could be 1.

\vfill\pagebreak

\section{The Details}  \label{syll sec}


\subsection{Details for case I}

Given the two equations
\begin{eqnarray}
vx&=&v'y\\
wz&=&w'y,
\end{eqnarray}
\textbf{reduce} the system to a single equation $R(x,y,z,v,v',w,w')=0$ where
\begin{equation}
R(x,y,z,v,v',w,w')\  :=\ (vx - v'y)^2 + (wz - w'y)^2.
\end{equation}
{\bf Eliminating} $y$ from $R=0$ yields $E(x,z,v,v',w,w') = 0$ where
\begin{eqnarray}
E(x,z,v,v',w,w') &:=& R(x,1,z,v,v',w,w')    R(x,0,z,v,v',w,w')\nonumber\\
&=&\Big((v   x-v')^2\ +\ (wz-w')^2\Big)   (v   x+wz).
\end{eqnarray}

The goal is to solve $E=0$ for $x$, for $1-x$ and for $vx$.

\subsection{Solving for $\pmb x$}\label{Case I x}

To solve $E=0$ for $x$, first express $E$ in the form $E_1x + E_0(1-x)$ by expanding it about $x$:
\begin{eqnarray*}
&&E(x,z,v,v',w,w')\\
&=& E(1,z,v,v',w,w')x +E(0,z,v,v',w,w')(1-x)\\
&=&\bigg(\Big((v-v')^2+(wz-w')^2\Big)   (v+wz)\bigg)x
\ +  \bigg(\Big(v'+(wz-w')^2\Big)   wz \bigg)(1-x).
\end{eqnarray*}
From $E_1x + E_0(1-x)=0$ one has $(E_0 - E_1)x = E_0$, and thus
one has the solution
\begin{eqnarray}{ }
  x & =& \frac{E(0,z,v,v',w,w') }{E(0,z,v,v',w,w')-E(1,z,v,v',w,w')}\nonumber\\
    & =& \frac{\Big(v'+(wz-w')^2\Big)   wz }{\Big(v'+(wz-w')^2\Big)   wz  -\Big((v-v')^2+(wz-w')^2\Big)   (v+wz )}  . \label{x-sol}
\end{eqnarray}

Now expand the right side of \eqref{x-sol} into a linear combination of constituents using
the following table---this is Boole's method of solution.\footnote
 { The coefficient $0/0$  represents an \textit{indefinite class}.
  $\infty$ is used as the value of any coefficient not in the form $0/b$ or $b/b$---Boole 
  preferred to use $1/0$ instead of the symbol $\infty$.

 The sum of the constituents with coefficient $\infty$  set equal to 0 gives the expansion
 of the constraint equation $E_1\cdot E_0 = 1$.}

\vfill\pagebreak

$$ \small
\begin{array}{l c c c c c r@{/}l c l}
&z &  v &  v' & w &  w' & \multicolumn{2}{c}{Coeff} & Value & \quad  Constituent\\
\hline
 1. & 0  &  0  & 0  & 0  & 0  &  0&0 &  0/0 &  (1-z)\cdot (1-v)\cdot (1-v')\cdot (1-w)\cdot (1-w') \\
 2. & 0 &  0  & 0  & 0  & 1  &  0&0 &  0/0 &  (1-z)\cdot (1-v)\cdot (1-v')\cdot (1-w)\cdot w' \\
 3. & 0 &  0  & 0  & 1  & 0  &  0&0 &  0/0 &  (1-z)\cdot (1-v)\cdot (1-v')\cdot w\cdot (1-w') \\
 4. & 0 &  0  & 0  & 1  & 1  &  0&0 &  0/0 &  (1-z)\cdot (1-v)\cdot (1-v')\cdot w\cdot w' \\
 5. & 0 &  0  & 1  & 0  & 0  &  0&0 &  0/0 &  (1-z)\cdot (1-v)\cdot v'\cdot (1-w)\cdot (1-w') \\
 6. & 0 &  0  & 1  & 0  & 1  &  0&0 &  0/0 &  (1-z)\cdot (1-v)\cdot v'\cdot (1-w)\cdot w' \\
 7. & 0 &  0  & 1  & 1  & 0  &  0&0 &  0/0 &  (1-z)\cdot (1-v)\cdot v'\cdot w\cdot (1-w') \\
 8. & 0 &  0  & 1  & 1  & 1  &  0&0 &  0/0 &  (1-z)\cdot (1-v)\cdot v'\cdot w\cdot w' \\
 9. & 0 &  1  & 0  & 0  & 0  &  0&1 &  0 &  (1-z)\cdot v\cdot (1-v')\cdot (1-w)\cdot (1-w') \\
10. & 0 &  1  & 0  & 0  & 1  &  0&2 &  0 &  (1-z)\cdot v\cdot (1-v')\cdot (1-w)\cdot w' \\
11. & 0 &  1  & 0  & 1  & 0  &  0&1 &  0 &  (1-z)\cdot v\cdot (1-v')\cdot w\cdot (1-w') \\
12. & 0 &  1  & 0  & 1  & 1  &  0&2 &  0 &  (1-z)\cdot v\cdot (1-v')\cdot w\cdot w' \\
13. & 0 &  1  & 1  & 0  & 0  &  0&0 &  0/0 &  (1-z)\cdot v\cdot v'\cdot (1-w)\cdot (1-w') \\
14. & 0 &  1  & 1  & 0  & 1  &  0&1 &  0 &  (1-z)\cdot v\cdot v'\cdot (1-w)\cdot w' \\
15. & 0 &  1  & 1  & 1  & 0  &  0&0 &  0/0 &  (1-z)\cdot v\cdot v'\cdot w\cdot (1-w') \\
16. & 0 &  1  & 1  & 1  & 1  &  0&1 &  0 &  (1-z)\cdot v\cdot v'\cdot w\cdot w' \\
17. & 1 &  0  & 0  & 0  & 0  &  0&0 &  0/0 &  z\cdot (1-v)\cdot (1-v')\cdot (1-w)\cdot (1-w') \\
18. & 1 &  0  & 0  & 0  & 1  &  0&0 &  0/0 &  z\cdot (1-v)\cdot (1-v')\cdot (1-w)\cdot w' \\
19. & 1 &  0  & 0  & 1  & 0  & -1&0 &  \infty &  z\cdot (1-v)\cdot (1-v')\cdot w\cdot (1-w') \\
20. & 1 &  0  & 0  & 1  & 1  &  0&0 &  0/0 &  z\cdot (1-v)\cdot (1-v')\cdot w\cdot w' \\
21. & 1 &  0  & 1  & 0  & 0  &  0&0 &  0/0 &  z\cdot (1-v)\cdot v'\cdot (1-w)\cdot (1-w') \\
22. & 1 &  0  & 1  & 0  & 1  &  0&0 &  0/0 &  z\cdot (1-v)\cdot v'\cdot (1-w)\cdot w' \\
23. & 1 &  0  & 1  & 1  & 0  & -2&0 &  \infty &  z\cdot (1-v)\cdot v'\cdot w\cdot (1-w') \\
24. & 1 &  0  & 1  & 1  & 1  & -1&0 &  \infty &  z\cdot (1-v)\cdot v'\cdot w\cdot w' \\
25. & 1 &  1  & 0  & 0  & 0  &  0&1 &  0 &  z\cdot v\cdot (1-v')\cdot (1-w)\cdot (1-w') \\
26. & 1 &  1  & 0  & 0  & 1  &  0&2 &  0 &  z\cdot v\cdot (1-v')\cdot (1-w)\cdot w' \\
27. & 1 &  1  & 0  & 1  & 0  & -1&3 &  \infty &  z\cdot v\cdot (1-v')\cdot w\cdot (1-w') \\
28. & 1 &  1  & 0  & 1  & 1  &  0&2 &  0 &  z\cdot v\cdot (1-v')\cdot w\cdot w' \\
29. & 1 &  1  & 1  & 0  & 0  &  0&0 &  0/0 &  z\cdot v\cdot v'\cdot (1-w)\cdot (1-w') \\
30. & 1 &  1  & 1  & 0  & 1  &  0&1 &  0 &  z\cdot v\cdot v'\cdot (1-w)\cdot w' \\
31. & 1 &  1  & 1  & 1  & 0  & -2&0 &  \infty &  z\cdot v\cdot v'\cdot w\cdot (1-w') \\
32. & 1 &  1  & 1  & 1  & 1  & -1&-1 & 1 &  z\cdot v\cdot v'\cdot w\cdot w' \\
\end{array}
$$
\centerline{Case I: Table for Equation \eqref{x-sol} \qquad}

\vfill\pagebreak

 The solution for $x$ is that it equals the sum of the \textit{coefficient $\times$ constituent} 
 for which the coefficient is 1 or 0/0. The same comment applies to the subsequent solutions
  for $1-x$ and $vx$.
  
Here is the expansion of \eqref{x-sol} in the form $f_1(v,v',w,w')z + f_0(v,v',w,w')(1-z)$, where
the numbers in parentheses (1),\ldots,(32) give the row numbers for the constituents
included in the expression for $x$. 
Next underbraces show which constituents give the 
indicated term. 
$x$, $z$ and $1-z$ are displayed in boldface to assist the reader in parsing
the expressions.
\begin{eqnarray}
\pmb{x}&=&\bigg( \big(32) + \frac{0}{0}\Big((17) + (18) + (20)+(21)+(22)+(29)\Big)\bigg) \pmb{z}   \nonumber\\
&+& 
\frac{0}{0}\Big((1) +   \cdots + (8) + (13) + (15)\Big)  \pmb{ (1-z)}\nonumber\\
&=&\bigg(\underbrace{v   v'   w   w'}_{(32)}  \ +\ 
 \frac{0}{0}\Big(\underbrace{(1-v)   (1-v')   w   w'}_{(20)}
\ +\  \underbrace{(1-v)  (1-w)}_{(17)+(18)+(21)+(22)}\ +\ 
\underbrace{v   v'   (1-w)   (1-w')}_{(29)}\Big)\bigg)   \pmb{ z} \nonumber\\
&+&\frac{0}{0}\Big(\underbrace{(1-v)}_{(1)+  \cdots+(8)}\ +\ \underbrace{v   v'  (1-w')}_{(13)+(15)}\Big)  \pmb{ (1-z) }  \nonumber\\
&=&\bigg(v   v'   w   w'  \ +\ 
 \frac{0}{0}  \Big((1-v)   (1-v')   w   w'
\ +\  (1-v)  (1-w)\ +\ 
v   v'   (1-w)   (1-w')\Big)\bigg)   \pmb{z}\nonumber\\
&+& \frac{0}{0} \Big((1-v)\ +\ v   v'  (1-w')\Big) \pmb{ (1-z)} .
\label{Ip233}  \\\nonumber
\end{eqnarray}
\fbox{Formula \eqref{Ip233} is Boole's (I.) on p.~233 of LT.}
\bigskip

\subsection{Solving for $\pmb{1-x}$}

From $E_1x + E_0(1-x)=0$ one has $(E_1 - E_0)(1-x) = E_1$, thus
\begin{eqnarray}{ }
  1-x 
  & =& \frac{E(1,z,v,v',w,w') }{E(1,z,v,v',w,w')-E(0,z,v,v',w,w')}\nonumber\\
  & =& \frac{\Big((v-v')^2+(wz-w')^2\Big)   (v+wz ) }{\Big((v-v')^2+(wz-w')^2\Big)   (v+wz )-\Big(v'+(wz-w')^2\Big)   wz }\label{1-x-sol}
\end{eqnarray}
so construct the table for \eqref{1-x-sol}:
\vfill\pagebreak

$$ \small
\begin{array}{l c c c c c r@{/}l c l}
&z &  v &  v' & w &  w' & \multicolumn{2}{c}{Coeff} & Value & \quad  Constituent\\
\hline
 1. & 0  & 0  & 0  & 0  & 0  &  0 & 0  &  0/0   &  (1-z)\cdot  (1-v)\cdot  (1-v')\cdot  (1-w)\cdot  (1-w') \\ 
 2. & 0  & 0  & 0  & 0  & 1  &  0 & 0  &  0/0   &  (1-z)\cdot  (1-v)\cdot  (1-v')\cdot  (1-w)\cdot  w'    \\ 
 3. & 0  & 0  & 0  & 1  & 0  &  0 & 0  &  0/0   &  (1-z)\cdot  (1-v)\cdot  (1-v')\cdot  w\cdot  (1-w')    \\ 
 4. & 0  & 0  & 0  & 1  & 1  &  0 & 0  &  0/0   &  (1-z)\cdot  (1-v)\cdot  (1-v')\cdot  w\cdot  w'        \\ 
 5. & 0  & 0  & 1  & 0  & 0  &  0 & 0  &  0/0   &  (1-z)\cdot  (1-v)\cdot  v'\cdot  (1-w)\cdot  (1-w')    \\ 
 6. & 0  & 0  & 1  & 0  & 1  &  0 & 0  & 0/0 &  (1-z)\cdot  (1-v)\cdot  v'\cdot  (1-w)\cdot  w'        \\ 
 7. & 0  & 0  & 1  & 1  & 0  &  0 & 0  &  0/0   &  (1-z)\cdot  (1-v)\cdot  v'\cdot  w\cdot  (1-w')        \\ 
 8. & 0  & 0  & 1  & 1  & 1  &  0 & 0  & 0/0 &  (1-z)\cdot  (1-v)\cdot  v'\cdot  w\cdot  w'            \\ 
 9. & 0  & 1  & 0  & 0  & 0  &  1 & 1  &  1     &  (1-z)\cdot  v\cdot  (1-v')\cdot  (1-w)\cdot  (1-w')    \\ 
10. & 0  & 1  & 0  & 0  & 1  &  2 & 2  &  1     &  (1-z)\cdot  v\cdot  (1-v')\cdot  (1-w)\cdot  w'        \\ 
11. & 0  & 1  & 0  & 1  & 0  &  1 & 1  &  1     &  (1-z)\cdot  v\cdot  (1-v')\cdot  w\cdot  (1-w')        \\ 
12. & 0  & 1  & 0  & 1  & 1  &  2 & 2  &  1     &  (1-z)\cdot  v\cdot  (1-v')\cdot  w\cdot  w'            \\ 
13. & 0  & 1  & 1  & 0  & 0  &  0 & 0  &  0/0   &  (1-z)\cdot  v\cdot  v'\cdot  (1-w)\cdot  (1-w')        \\ 
14. & 0  & 1  & 1  & 0  & 1  &  1 & 1 &  1     &  (1-z)\cdot  v\cdot  v'\cdot  (1-w)\cdot  w'            \\ 
15. & 0  & 1  & 1  & 1  & 0  &  0 & 0  &  0/0   &  (1-z)\cdot  v\cdot  v'\cdot  w\cdot  (1-w')            \\ 
16. & 0  & 1  & 1  & 1  & 1  &  1 & 1 &  1     &  (1-z)\cdot  v\cdot  v'\cdot  w\cdot  w'                \\ 
17. & 1  & 0  & 0  & 0  & 0  &  0 & 0  &  0/0   &  z\cdot  (1-v)\cdot  (1-v')\cdot  (1-w)\cdot  (1-w')    \\ 
18. & 1  & 0  & 0  & 0  & 1  &  0 & 0  &  0/0   &  z\cdot  (1-v)\cdot  (1-v')\cdot  (1-w)\cdot  w'        \\ 
19. & 1  & 0  & 0  & 1  & 0  &  1 & 0  & \infty &  z\cdot  (1-v)\cdot  (1-v')\cdot  w\cdot  (1-w')        \\ 
20. & 1  & 0  & 0  & 1  & 1  &  0 & 0  &  0/0   &  z\cdot  (1-v)\cdot  (1-v')\cdot  w\cdot  w'            \\ 
21. & 1  & 0  & 1  & 0  & 0  &  0 & 0  &  0/0   &  z\cdot  (1-v)\cdot  v'\cdot  (1-w)\cdot  (1-w')        \\ 
22. & 1  & 0  & 1  & 0  & 1  &  0 & 0  & 0/0 &  z\cdot  (1-v)\cdot  v'\cdot  (1-w)\cdot  w'            \\ 
23. & 1  & 0  & 1  & 1  & 0  &  2 & 0  & \infty &  z\cdot  (1-v)\cdot  v'\cdot  w\cdot  (1-w')            \\ 
24. & 1  & 0  & 1  & 1  & 1  &  1& 0  &  \infty   &  z\cdot  (1-v)\cdot  v'\cdot  w\cdot  w'                \\ 
25. & 1  & 1  & 0  & 0  & 0  &  1 & 1  &  1     &  z\cdot  v\cdot  (1-v')\cdot  (1-w)\cdot  (1-w')        \\ 
26. & 1  & 1  & 0  & 0  & 1  &  2 & 2  &  1     &  z\cdot  v\cdot  (1-v')\cdot  (1-w)\cdot  w'            \\ 
27. & 1  & 1  & 0  & 1  & 0  &  4 & 3  & \infty &  z\cdot  v\cdot  (1-v')\cdot  w\cdot  (1-w')            \\ 
28. & 1  & 1  & 0  & 1  & 1  &  2 & 2  &  1     &  z\cdot  v\cdot  (1-v')\cdot  w\cdot  w'                \\ 
29. & 1  & 1  & 1  & 0  & 0  &  0 & 0  &  0/0   &  z\cdot  v\cdot  v'\cdot  (1-w)\cdot  (1-w')            \\ 
30. & 1  & 1  & 1  & 0  & 1  &  1 & 1 & 1     &  z\cdot  v\cdot  v'\cdot  (1-w)\cdot  w'                \\ 
31. & 1  & 1  & 1  & 1  & 0  &  2 & 0  & \infty &  z\cdot  v\cdot  v'\cdot  w\cdot  (1-w')                \\ 
32. & 1  & 1  & 1  & 1  & 1  &  0 & -1  & 0     &  z\cdot  v\cdot  v'\cdot  w\cdot  w'                    \\ 
\end{array}
$$

\centerline{Case I: Table for Equation \eqref{1-x-sol} \qquad}
\vfill\pagebreak

Abbreviating the term $v'(1-w) + (1-v')w$ by $v'\triangle w$, the table gives
\begin{eqnarray}
\pmb{1-x}&=&\bigg[(25)+(26) +(28)+  {(30)} + \frac{0}{0}\Big((17) + (18) + (20)+(21)+(22)+(29)\Big)\bigg]   \pmb{z} \nonumber\\
&+&
\bigg[(9) + (10) +(11) +(12) + (14) +(16)\nonumber\\
 &&\qquad +\ 
\frac{0}{0}   \Big(
(1) + (2) + (3)+(4)+(5)+ (7) + (13) + (15)\Big)\bigg]   \pmb{(1-z)}\nonumber\\
&=&\bigg[\underbrace{v   (1-v')   (1-w)}_{(25)+(26)}  \ +\ \underbrace{v   w'  (v'\triangle w) }_{(28)+(30)}\nonumber\\  
\nonumber\\
&&\qquad\qquad + \ 
\frac{0}{0}  \Big(  \underbrace{(1-v)  (1-w)}_{(17)+(18)+(21)+(22)}\ 
\underbrace{(1-v)   (1-v')   w    w'}_{(20)} \ +\  
\underbrace{v   v'  (1-w)  (1-w')}_{(29)} \Big)\bigg]   \pmb{z}\nonumber\\
&+& \bigg[ {\underbrace{v   (1-v')}_{(9)+\cdots+(12)}} \ +\   \underbrace{ v   {v'}   w'}_{(14)+(16)}\ +\ \frac{0}{0}  \Big( \underbrace{(1-v)}_{(1)+  \cdots+(8)} + 
\underbrace{v  v' (1-w')}_{(13)+(15)}\Big)   \bigg] \pmb{(1-z)}\nonumber\\
&=&\bigg[ v   (1-v')   (1-w)  \ +\  v   w'  (v'\triangle w) 
\nonumber\\
&& \qquad\qquad +\  \frac{0}{0}  \Big(  (1-v)  (1-w) +\ 
 (1-v)   (1-v')   w    w' \ +\  
 v   v'  (1-w)  (1-w') \Big)\bigg]   \pmb{z}\nonumber\\
&+& \bigg[ \RR{\pmb{ v   (1-v')}} \ +\     v \;  \RR{\pmb{v'}}   w'\ +\ \frac{0}{0}  \Big(  (1-v) \ +\  
 v   v'  (1-w')\Big)   \bigg]  \pmb{(1-z)} .
\label{IIp233}  \\\nonumber
\end{eqnarray}
\fbox{Formula \eqref{IIp233} is Boole's (II.) on p.~233 of LT} except for the red colored items: the term $v  (1-v')$ is completely missing, and in the next term the $v'$ in LT is given as $(1-w)$. \\

\subsection{Solving for $\pmb {vx}$}\label{Case I vx}

Multiplying \eqref{Ip233} by $v$ gives
\begin{eqnarray}
\pmb{vx}
&=&\Big(v   {v'}   w   w'  \ +\ 
 \frac{0}{0}\;v   v'   (1-w)   (1-w')\Big)   \pmb{z} \ 
 +\ \frac{0}{0}\; v \;  \RR{\pmb{v'}}  (1-w')   \pmb{(1-z)}.\label{IIIp233}   \\\nonumber
\end{eqnarray}
\fbox{Formula \eqref{IIIp233} is Boole's (III.) on p.~233 of LT} except that Boole omitted the
 $\RR{\pmb{v'}}$ indicated in red.
 
\vfill

\pagebreak

\subsection{Summary of Solutions for Case I}

\begin{eqnarray}
\pmb{x}&=&\bigg[v   v'   w   w'  \ +\ 
 \frac{0}{0}  \Big((1-v)   (1-v')   w   w'
\ +\  (1-v)  (1-w)\ +\ 
v   v'   (1-w)   (1-w')\Big)\bigg]   \pmb{z} \nonumber\\
&+& \frac{0}{0} \Big((1-v)\ +\ v   v'  (1-w')\Big)  \pmb{(1-z)}\label{Sol1} \\
\pmb{1-x}&=&\bigg[ v   (1-v')   (1-w)  \ +\  v   w'  (v'\triangle w) 
\nonumber\\
&& \qquad\qquad +\  \frac{0}{0}  \Big(  (1-v)  (1-w) +\ 
 (1-v)   (1-v')   w    w' \ +\  
 v   v'  (1-w)  (1-w') \Big)\bigg]   \pmb{z}\nonumber\\
&+& \bigg[ \pmb{\RR{ v   (1-v')}} \ +\     v\;   \RR{\pmb{v'}}   w'\ +\ \frac{0}{0}  \Big(  (1-v) + 
 v   v'  (1-w')\Big)   \bigg]   \pmb{(1-z)} .\label{Sol2} \\
 \pmb{vx} &=&\bigg(v   {v'}   w   w'  \ +\ 
 \frac{0}{0}\Big(v   v'   (1-w)   (1-w')\Big)\bigg)    \pmb{z}
\ +\ \frac{0}{0}\Big( v  \; \RR{\pmb{v'}}  (1-w')\Big)   \pmb{ (1-z)}. \label{Sol3}
\end{eqnarray}

We consider each of the equations \eqref{Sol1}--\eqref{Sol3} above in turn---the errors
 in LT noted in
red type do not affect the conclusions Boole derived regarding valid syllogisms.
\begin{itemize}
\item
For \eqref{Sol1}  the coefficient of $z$ cannot be made to vanish using permissible substitutions
of 1 for $v,v',w,w'$.
To make the coefficient of $1-z$ vanish one needs to assign $v=w'=1$. 
Then \eqref{Sol1}  reduces to $x = v'wz$.\\

\item
For \eqref{Sol2}  the only possibility is $v'=w=1$ since one cannot force the coefficient of 
$1-z$ to be 0 with a permissible substitution.
 Then \eqref{Sol2} reduces to 
$$1-x =
\Big(v w'\ +\ \dfrac{0}{0} (1-vw')  \Big)  (1-z). $$
\smallskip

\item
For \eqref{Sol3}  one must set $w'=1$, and then \eqref{Sol3} reduces $vx = vv'wz$.
One can additionally set $v'=1$, reducing \eqref{Sol3} to $vx = vwz$.

\end{itemize}
\medskip

\vfill\pagebreak

\subsection{Details for Case II }

Given the two equations
\begin{eqnarray}
vx&=&v'y\\
wz&=&w'(1-y)
\end{eqnarray}
\textbf{reduce} the system to a single equation $R(x,y,z,v,v',w,w')=0$ where
\begin{equation}
R(x,y,z,v,v',w,w')\  :=\ (vx - v'y)^2 + \Big(wz - w'(1-y)\Big)^2.
\end{equation}
{\bf Eliminating} $y$ from $R=0$ yields $E(x,z,v,v',w,w') = 0$ where
\begin{eqnarray}
E(x,z,v,v',w,w') &:=& R(x,1,z,v,v',w,w')    R(x,0,z,v,v',w,w')\nonumber\\
&=&\Big((vx -v')^2\ +\  wz\Big)   \Big(vx + (wz -w')^2\Big).
\end{eqnarray}

Now we want to solve $E=0$ for $x$, for $1-x$ and for $vx$.

\subsection{Solving for $x$}\label{Case II x}

To solve $E=0$ for $x$, first express $E$ in the form $E_1x + E_0(1-x)$ by expanding it about $x$:
\begin{eqnarray*}
&&E(x,z,v,v',w,w')\\
&=& E(1,z,v,v',w,w')x +E(0,z,v,v',w,w')(1-x)\\
&=&\Big((v-v')^2 + w z\Big)   \Big(v+(wz -w')^2\Big)x
\ + \ (v'+ wz) (wz-w')^2(1-x).
\end{eqnarray*}
From $E_1x + E_0(1-x)=0$ one has $(E_0-E_1)x =E_0$, thus
\begin{eqnarray}{ }
  x & =& \frac{E(0,z,v,v',w,w') }{E(0,z,v,v',w,w')-E(1,z,v,v',w,w')}\nonumber\\
    & =& \frac{(v'+ wz) (wz-w')^2}{
    (v'+ wz) (wz-w')^2 - \Big((v-v')^2 + wz\Big) \Big(v+(wz -w')^2\Big) } .
    \label{Case II x-sol}
\end{eqnarray}

Constructing the table for the right side of \eqref{Case II x-sol} gives:
\vfill \pagebreak

$$ \small
\begin{array}{l c c c c c r@{/}l c l}
&z &  v &  v' & w &  w' & \multicolumn{2}{c}{Coeff} & Value & \quad  Constituent\\
\hline
 1. & 0  & 0  & 0  & 0  & 0  &  0&0  &  0/0   &  (1-z)(1-v)(1-v')(1-w)(1-w') \\ 
 2. & 0  & 0  & 0  & 0  & 1  &  0&0  &  0/0   &  (1-z)(1-v)(1-v')(1-w)w'    \\ 
 3. & 0  & 0  & 0  & 1  & 0  &  0&0  &  0/0   &  (1-z)(1-v)(1-v')w(1-w')    \\ 
 4. & 0  & 0  & 0  & 1  & 1  &  0&0  &  0/0   &  (1-z)(1-v)(1-v')ww'        \\ 
 5. & 0  & 0  & 1  & 0  & 0  &  0&0  &  0/0   &  (1-z)(1-v)v'(1-w)(1-w')    \\ 
 6. & 0  & 0  & 1  & 0  & 1  & -1&0  & \infty &  (1-z)(1-v)v'(1-w)w'        \\ 
 7. & 0  & 0  & 1  & 1  & 0  &  0&0  &  0/0   &  (1-z)(1-v)v'w(1-w')        \\ 
 8. & 0  & 0  & 1  & 1  & 1  & -1&0  & \infty &  (1-z)(1-v)v'ww'            \\ 
 9. & 0  & 1  & 0  & 0  & 0  &  0&1  &  0     &  (1-z)v(1-v')(1-w)(1-w')    \\ 
10. & 0  & 1  & 0  & 0  & 1  &  0&2  &  0     &  (1-z)v(1-v')(1-w)w'        \\ 
11. & 0  & 1  & 0  & 1  & 0  &  0&1  &  0     &  (1-z)v(1-v')w(1-w')        \\ 
12. & 0  & 1  & 0  & 1  & 1  &  0&2  &  0     &  (1-z)v(1-v')ww'            \\ 
13. & 0  & 1  & 1  & 0  & 0  &  0&0  &  0/0   &  (1-z)vv'(1-w)(1-w')        \\ 
14. & 0  & 1  & 1  & 0  & 1  & -1&-1 &  1     &  (1-z)vv'(1-w)w'            \\ 
15. & 0  & 1  & 1  & 1  & 0  &  0&0  &  0/0   &  (1-z)vv'w(1-w')            \\ 
16. & 0  & 1  & 1  & 1  & 1  & -1&-1 &  1     &  (1-z)vv'ww'                \\ 
17. & 1  & 0  & 0  & 0  & 0  &  0&0  &  0/0   &  z(1-v)(1-v')(1-w)(1-w')    \\ 
18. & 1  & 0  & 0  & 0  & 1  &  0&0  &  0/0   &  z(1-v)(1-v')(1-w)w'        \\ 
19. & 1  & 0  & 0  & 1  & 0  & -1&0  & \infty &  z(1-v)(1-v')w(1-w')        \\ 
20. & 1  & 0  & 0  & 1  & 1  &  0&0  &  0/0   &  z(1-v)(1-v')ww'            \\ 
21. & 1  & 0  & 1  & 0  & 0  &  0&0  &  0/0   &  z(1-v)v'(1-w)(1-w')        \\ 
22. & 1  & 0  & 1  & 0  & 1  & -1&0  & \infty &  z(1-v)v'(1-w)w'            \\ 
23. & 1  & 0  & 1  & 1  & 0  & -2&0  & \infty &  z(1-v)v'w(1-w')            \\ 
24. & 1  & 0  & 1  & 1  & 1  &  0&0  &  0/0   &  z(1-v)v'ww'                \\ 
25. & 1  & 1  & 0  & 0  & 0  &  0&1  &  0     &  zv(1-v')(1-w)(1-w')        \\ 
26. & 1  & 1  & 0  & 0  & 1  &  0&2  &  0     &  zv(1-v')(1-w)w'            \\ 
27. & 1  & 1  & 0  & 1  & 0  & -1&3  & \infty &  zv(1-v')w(1-w')            \\ 
28. & 1  & 1  & 0  & 1  & 1  &  0&2  &  0     &  zv(1-v')ww'                \\ 
29. & 1  & 1  & 1  & 0  & 0  &  0&0  &  0/0   &  zvv'(1-w)(1-w')            \\ 
30. & 1  & 1  & 1  & 0  & 1  & -1&-1 &  1     &  zvv'(1-w)w'                \\ 
31. & 1  & 1  & 1  & 1  & 0  & -2&0  & \infty &  zvv'w(1-w')                \\ 
32. & 1  & 1  & 1  & 1  & 1  &  0&1  &  0     &  zvv'ww'  .
\end{array}
$$
\vspace*{0.2 em}

\centerline{Case II: Table for Equation \eqref{Case II x-sol} \qquad}
\vfill\pagebreak

Thus
\begin{eqnarray}
 \pmb{x}&=&
\bigg(\ (30)
 \ +\  
\frac{0}{0}
\Big(
(17) + (18) + (20)+(21)+(24)+(29)
\Big)
\bigg)  \pmb{z} \nonumber\\
&+&\bigg(
 (14) + (16)  
 \ +\ 
\frac{0}{0}
\Big( 
 (1)+\cdots+(5) + (7) + (13) + (15)
 \Big) 
  \pmb{(1-z)} \nonumber\\
&=&\bigg(
\underbrace{v   v'   (1-w)   w'}_{(30)}  \ +\ 
 \frac{0}{0}
 \Big(\underbrace{(1-v)   (1-v')(1-w)}_{(17)+(18)}
\ +\ 
\underbrace{(1-v)ww'}_{(20)+(24)}\ +\
\underbrace{v'   (1-w) (1-w')}_{(21)+(29)}\Big)\bigg)    \pmb{z} \nonumber\\
&+&  
\bigg(
 \underbrace{vv'w'}_{(14)+(16)}
\ +\frac{0}{0}
\Big( 
\underbrace{(1-v)(1-v')}_{(1)+  \cdots + (4)} \ +\ 
\underbrace{v'  (1-w')}_{(5)+(7) +(13)+(15)}
\Big)  
\bigg) 
 \pmb{(1-z)}   \nonumber\\
&&\nonumber\\&&\nonumber\\
&=&\bigg(
{v   v'   (1-w)   w'}  \ +\ 
 \frac{0}{0}
 \Big(
{(1-v)   (1-v')(1-w)}
\ +\ 
{(1-v)ww'}\ +\
{v'   (1-w) (1-w')}
\Big)\bigg)   \pmb{z} \nonumber\\
&+& 
\bigg(
 vv'w'
\ +\frac{0}{0}
\Big( 
{(1-v)(1-v')} \ +\ 
{v'  (1-w')}
\Big)  
\bigg) 
 \pmb{(1-z). }
\label{ICp234}   \\\nonumber 
\end{eqnarray}
\fbox{Formula \eqref{ICp234} is Boole's (IV.) on p.~234 of LT.}
\bigskip

\subsection{Solving for $1-x$}

From $(E_1 - E_0)(1-x) = E_1$ one has
\begin{eqnarray}{ }
  1-x & =& \frac{E(1,z,v,v',w,w') }{E(1,z,v,v',w,w')-E(0,z,v,v',w,w')}   \nonumber\\
    & =& \frac{ \Big((v-v')^2 + w z\Big)   \Big(v+(wz -w')^2\Big)}{
    \Big((v-v')^2 + w z\Big)   \Big(v+(wz -w')^2\Big) - (v'+ wz) (wz-w')^2}
    \label{Case II 1-x-sol}
\end{eqnarray}
so construct the table for \eqref{Case II 1-x-sol}:
\vfill
\pagebreak

$$\small
\begin{array}{l c c c c c r@{/}l c l}
&z &  v &  v' & w &  w' & \multicolumn{2}{c}{Coeff} & Value & \quad  Constituent\\
\hline  
 1. & 0  & 0  & 0  & 0  & 0  &  0&0  &  0/0   &  (1-z)(1-v)(1-v')(1-w)(1-w) \\ 
 2. & 0  & 0  & 0  & 0  & 1  &  0&0  &  0/0   &  (1-z)(1-v)(1-v')(1-w)w      \\ 
 3. & 0  & 0  & 0  & 1  & 0  &  0&0  &  0/0   &  (1-z)(1-v)(1-v')w(1-w)      \\ 
 4. & 0  & 0  & 0  & 1  & 1  &  0&0  &  0/0   &  (1-z)(1-v)(1-v')ww            \\ 
 5. & 0  & 0  & 1  & 0  & 0  &  0&0  &  0/0   &  (1-z)(1-v)v'(1-w)(1-w)      \\ 
 6. & 0  & 0  & 1  & 0  & 1  &  1&0  & \infty &  (1-z)(1-v)v'(1-w)w          \\ 
 7. & 0  & 0  & 1  & 1  & 0  &  0&0  &  0/0   &  (1-z)(1-v)v'w(1-w)            \\ 
 8. & 0  & 0  & 1  & 1  & 1  &  1&0  & \infty &  (1-z)(1-v)v'ww              \\ 
 9. & 0  & 1  & 0  & 0  & 0  &  1&1  &  1     &  (1-z)v(1-v')(1-w)(1-w)    \\ 
10. & 0  & 1  & 0  & 0  & 1  &  2&2  &  1     &  (1-z)v(1-v')(1-w)w        \\ 
11. & 0  & 1  & 0  & 1  & 0  &  1&1  &  1     &  (1-z)v(1-v')w(1-w)        \\ 
12. & 0  & 1  & 0  & 1  & 1  &  2&2  &  1     &  (1-z)v(1-v')ww             \\ 
13. & 0  & 1  & 1  & 0  & 0  &  0&0  &  0/0   &  (1-z)vv'(1-w)(1-w)        \\ 
14. & 0  & 1  & 1  & 0  & 1  &  0&-1 &  0     &  (1-z)vv'(1-w)w            \\ 
15. & 0  & 1  & 1  & 1  & 0  &  0&0  &  0/0   &  (1-z)vv'w(1-w)              \\ 
16. & 0  & 1  & 1  & 1  & 1  &  0&-1 &  0     &  (1-z)vv'ww                  \\ 
17. & 1  & 0  & 0  & 0  & 0  &  0&0  &  0/0   &  z(1-v)(1-v')(1-w)(1-w')    \\ 
18. & 1  & 0  & 0  & 0  & 1  &  0&0  &  0/0   &  z(1-v)(1-v')(1-w)w'        \\ 
19. & 1  & 0  & 0  & 1  & 0  &  1&0  & \infty &  z(1-v)(1-v')w(1-w')        \\ 
20. & 1  & 0  & 0  & 1  & 1  &  0&0  &  0/0   &  z(1-v)(1-v')ww'            \\ 
21. & 1  & 0  & 1  & 0  & 0  &  0&0  &  0/0   &  z(1-v)v'(1-w)(1-w')        \\ 
22. & 1  & 0  & 1  & 0  & 1  &  1&0  & \infty &  z(1-v)v'(1-w)w'            \\ 
23. & 1  & 0  & 1  & 1  & 0  &  2&0  & \infty &  z(1-v)v'w(1-w')            \\ 
24. & 1  & 0  & 1  & 1  & 1  &  0&0  &  0/0   &  z(1-v)v'ww'                \\ 
25. & 1  & 1  & 0  & 0  & 0  &  1&1  &  1     &  zv(1-v')(1-w)(1-w')        \\ 
26. & 1  & 1  & 0  & 0  & 1  &  2&2  &  1     &  zv(1-v')(1-w)w'            \\ 
27. & 1  & 1  & 0  & 1  & 0  &  4&3  & \infty &  zv(1-v')w(1-w')            \\ 
28. & 1  & 1  & 0  & 1  & 1  &  2&2  &  1     &  zv(1-v')ww'                \\ 
29. & 1  & 1  & 1  & 0  & 0  &  0&0  &  0/0   &  zvv'(1-w)(1-w')            \\ 
30. & 1  & 1  & 1  & 0  & 1  &  0&-1 &  0     &  zvv'(1-w)w'                \\ 
31. & 1  & 1  & 1  & 1  & 0  &  2&0  & \infty &  zvv'w(1-w')                \\ 
32. & 1  & 1  & 1  & 1  & 1  &  1&1  &  1     &  zvv'ww'                    \\ 

\end{array}
$$
\vspace*{0.2 em}

\centerline{Case II: Table for Equation \eqref{Case II 1-x-sol} \qquad}
\vfill\pagebreak

The table gives
\begin{eqnarray}
 \pmb{1-x}&=&\bigg(\ (25) + (26)+(28)+(32)
 \ +\  
\frac{0}{0}
\Big(
(17) + (18) + (20)+(21)+(24)+(29)
\Big)
\bigg)  \pmb{z} \nonumber\\
&+&\bigg(
(9)+\cdots + (12) 
 \ +\ 
\frac{0}{0}
\Big( 
 (1)+\cdots+(5) + (7) + (13) + (15)
 \Big) 
  \pmb{(1-z)} \nonumber\\
&=&\bigg(
\underbrace{v   (1-v')   (1-w) }_{(25) + (26)}  \ +\ 
\underbrace{v   w  w' }_{(28)+(32)} \nonumber\\
&&\quad  \ +\ 
 \frac{0}{0}
 \Big(\underbrace{(1-v)   (1-v')(1-w)}_{(17)+(18)}
\ +\ 
\underbrace{(1-v)ww'}_{(20)+(24)}\ +\
\underbrace{v'   (1-w) (1-w')}_{(21)+(29)}\Big)\bigg)    \pmb{z}  \nonumber\\
&+&
\bigg(
 \underbrace{v(1-v')}_{(9)+\cdots + (12)}
\ +\frac{0}{0}
\Big( 
\underbrace{(1-v)(1-v')}_{(1)+  \cdots + (4)} \ +\ 
\underbrace{v'  (1-w')}_{(5)+(7) +(13)+(15)}
\Big)  
\bigg) 
 \pmb{(1-z)}   \nonumber\\
&&\nonumber\\
&=&\bigg(
{v   (1-v')   (1-w) }  \ +\ 
{v   w  w' } \nonumber\\
&&\quad  \ +\ 
 \frac{0}{0}
 \Big(
 {(1-v)   (1-v')(1-w)}
\ +\ 
{(1-v)ww'}\ +\
{v'   (1-w) (1-w')}\Big)\bigg)    \pmb{z} \nonumber\\
&+& 
\bigg(
{v(1-v')}
\ +\frac{0}{0}
\Big( 
{(1-v)(1-v')} \ +\ 
{v'  (1-w')}
\Big)  
\bigg) 
 \pmb{(1-z)}.  \label{IICp235}     \\\nonumber
\end{eqnarray}
\fbox{Formula \eqref{IICp235} is Boole's (V.) on p.~235 of LT.} 
\bigskip

\subsection{Solving for $vx$}\label{Case II vx}

Multiplying \eqref{ICp234} by $v$ gives
\begin{eqnarray}
 \pmb{vx}
&=&
\Big(
{v   v'   (1-w)   w'}  \ +\ 
 \frac{0}{0}
{vv'   (1-w) (1-w')}
\Big)    \pmb{z} \nonumber\\
&+&  
\Big(
 vv'w'
\ +\frac{0}{0} 
{vv'  (1-w')}
\Big)  
 \pmb{(1-z) }.
 \label{IIICp235}  \\\nonumber
\end{eqnarray}
\fbox{Formula \eqref{IIICp235} is Boole's (VI.) on p.~235 of LT.} 

\vfill\pagebreak

\subsection{Summary of Solutions for Case II}

\begin{eqnarray}
 \pmb{x}&=&\bigg(
{v   v'   (1-w)   w'}  \ +\ 
 \frac{0}{0}
 \Big(
{(1-v)   (1-v')(1-w)}
\ +\ 
{(1-v)ww'}\ +\
{v'   (1-w) (1-w')}
\Big)\bigg)    \pmb{z} \nonumber\\
&+&
\bigg(
 vv'w'
\ +\frac{0}{0}
\Big( 
{(1-v)(1-v')} \ +\ 
{v'  (1-w')}
\Big)  
\bigg) 
 \pmb{(1-z)}
\label{Sol1C} \\
 \pmb{1-x}&=&\bigg(
{v   (1-v')   (1-w) }  \ +\ 
{v   w  w' } \nonumber\\
&&\qquad  \ +\ 
 \frac{0}{0}
 \Big(
 {(1-v)   (1-v')(1-w)}
\ +\ 
{(1-v)ww'}\ +\
{v'   (1-w) (1-w')}\Big)\bigg)    \pmb{z} \nonumber\\
&+&
\bigg(
{v(1-v')}
\ +\frac{0}{0}
\Big( 
{(1-v)(1-v')} \ +\ 
{v'  (1-w')}
\Big)  
\bigg) 
 \pmb{(1-z) }.
\label{Sol2C} \\
 \pmb{vx}
&=&
\Big(
{v   v'   (1-w)   w'}  \ +\ 
 \frac{0}{0}
{vv'   (1-w) (1-w')}
\Big)    \pmb{z} \nonumber\\
&+& 
\Big(
 vv'w'
\ +\frac{0}{0} 
{vv'  (1-w')}
\Big)  
 \pmb{(1-z)}.
\label{Sol3C}
\end{eqnarray}

We consider each of the equations \eqref{Sol1C}--\eqref{Sol3C} above in turn.\\

\begin{itemize}
\item
For \eqref{Sol1C}  the coefficient of $1-z$ cannot be made to vanish 
using a permissible substitution. 
To make the coefficient of $z$ vanish one needs to assign $v=w=1$. 
Then \eqref{Sol1C}  reduces to $x =\Big(v'w'+ \dfrac{0}{0}v'(1-w')\Big)(1-z)$. \\
\item
For \eqref{Sol2C}  the coefficient of $z$ cannot be made to vanish using a permissible substitution. 
To make the coefficient of $1-z$ vanish one needs to assign $v'=w'=1$. 
Then \eqref{Sol2C}  reduces to $1-x =\Big(vw+ \dfrac{0}{0}(1-v)w\Big)z$.\\
\item
For \eqref{Sol3C}  the coefficient of $1-z$ cannot be made to vanish using a permissible substitution. To make the coefficient of $z$ vanish one needs to assign $w=1$. 
Then \eqref{Sol3C}  reduces to  $vx =\Big(vv'w'+ \dfrac{0}{0}vv'(1-w')\Big)(1-z)$.
One can additionally set $v'=1$, giving the reduction to $vx =\Big(vw'+ \dfrac{0}{0}v(1-w')\Big)(1-z)$.\\\end{itemize}
\vfill \pagebreak

\section{Relevant Equational Arguments}

The above justifies the following equational arguments: \\

For Case I:

\noindent
$
\begin{array}{r c l}
{\bf (a)}\quad x &=& vy\\
wz &=& y\\
\therefore x &=& vwz
\end{array}
\qquad
\begin{array}{r c l}
{\bf (b)}\quad vx &=& y\\
z &=& wy\\
\therefore 1-x &=& \big(v w\ +\ \frac{0}{0} (1-vw)  \big)  (1-z) \\
&&
\end{array} \\ 
\begin{array}{r c l}
{\bf (c)}\quad vx &=& vy\\
wz &=& y\\
\therefore vx &=& vwz
\end{array}
\qquad
\begin{array}{r c l}
{\bf (d)}\quad vx &=& y\\
wz &=& y\\
\therefore vx &=& vwz .
\end{array}
$
\bigskip

From Case I  Boole claimed for categorical premises with like middle terms (p.~234 of LT):\footnote
{These three rules from Cases I and II, for determining valid syllogisms,
were originally announced as
the main contribution of Boole's 1848 paper \cite{Boole-1848} {\em The Calculus of Logic}.
Boole was quite pleased that he had abolished the need for the Aristotelian concepts of 
{\em figure} and {\em mood} after redefining the notions of {\em quality} and {\em quantity}
to apply to terms instead of to propositions. } 

\begin{quote}
\textsc{Condition of Inference.}--- {One middle term, at least, universal.}\\
\textsc{Rule of Inference.}--- {Equate the extremes.}
\end{quote}
\bigskip

For Case II:

$
\begin{array}{r c l}
{\bf (a)}\quad x &=& vy\\
z &=& w(1-y)\\
\therefore x &=&\big(vw+ \frac{0}{0}v(1-w)\big)(1-z)
\end{array}
\qquad
\begin{array}{r c l}
{\bf (b)}\quad vx &=& y\\
wz &=& 1-y\\
\therefore 1-x &=& \big(vw+ \frac{0}{0}(1-v)w\big)z\\
&&
\end{array}\\
\begin{array}{r c l}
{\bf (c)}\quad vx &=& vy\\
z &=& w(1-y)\\
\therefore vx &=&\big(vw+ \frac{0}{0}v(1-w)\big)(1-z)
\end{array}
\qquad
\begin{array}{r c l}
{\bf (d)}\quad vx &=& y\\
z &=& w(1-y)\\
\therefore vx &=& \big(vw+ \frac{0}{0}v(1-w)\big)(1-z).
\end{array}
$ \bigskip

From Case II Boole claimed (pp.~235, 236 of LT) for premises with unlike 
middle terms:\footnote{Boole mistakenly said the first of 
these two rules applied when the middle terms had \textit{like} quality.}
\medskip
\begin{quote}
\textsc{First Condition of Inference.}--- \textit{At least one universal extreme.} \\
\textsc{Rule of Inference.}--- \textit{Change the quantity and quality of that extreme,
and equate the result to the other extreme.}
\end{quote}
\medskip
\begin{quote}
\textsc{Second Condition of Inference.}--- \textit{Two universal middle terms.} \\
\textsc{Rule of Inference.}--- \textit{Change the quantity and quality of either extreme,
and equate the result to the other extreme unchanged.} \\
\end{quote}

\section{The Advantage of Using Negated Equations}

Note that all of Boole's categorical propositions $\Phi(X,Y)$ can be expressed in either 
the form $\alpha\beta = 0$ or $\alpha\beta \neq 0$, where $\alpha$ is either $x$
or $1-x$ and $\beta$ is either $y$ or $1-y$. For example, All not-$X$ is $Y$ is expressed
by $(1-x)(1-y) =0$, and Some not-$X$ is $Y$ is expressed by $(1-x)y \neq0$.

By using negated equations as well as equations the mathematical description of valid
syllogisms can be condensed into four cases where $\overline{\alpha}$ changes $x$ to
$1-x$ and vice-versa, etc.
$$
\begin{array}{r c l}
\alpha \gamma &=& 0\\
\beta \gamma &=&0\\
\hline
\therefore \ 
\overline{\alpha} \overline{\beta} &\neq& 0
\end{array}
\qquad
\begin{array}{r c l}
\alpha \gamma &=& 0\\
\beta \overline{\gamma} &=&0\\
\hline
\therefore \ 
{\alpha} {\beta} &=& 0
\end{array}
\qquad
\begin{array}{r c l}
\alpha \gamma &=& 0\\
\beta \gamma &\neq&0\\
\hline
\therefore \ 
\overline{\alpha} {\beta} &\neq& 0
\end{array}
\qquad
\begin{array}{r c l}
\alpha \gamma &=& 0\\
\beta \overline{\gamma} &\neq&0\\
\hline
\therefore \ 
{\alpha} {\beta} &\neq& 0
\end{array}
$$
These arguments are easily confirmed using Venn diagrams, just as the two forms of premises
$$
\begin{array}{r c l}
\alpha \gamma &\neq& 0\\
\beta \gamma &\neq&0\\
\end{array}
\qquad \qquad
\begin{array}{r c l}
\alpha \gamma &\neq& 0\\
\beta \overline{\gamma} &\neq&0\\
\end{array}
$$
can be shown to have no valid conclusion. (A formal proof theory for dealing with $\neq$, however, 
requires more work.)

Boole limited his general theory to universal propositions for a good reason---there 
was no  general elimination theorem if
one included particular propositions. For example consider the premises\smallskip

Some $Y$ is $X$\\
\indent
Some not-$Y$ is $X$.\smallskip

\noindent
The result of eliminating $Y$ is that $X$ has at least two elements. However this 
conclusion cannot
be expressed in his algebra of logic.\\

{\sc References}   \label{references}

\end{document}